\documentclass[11pt]{amsart}
\usepackage{amsmath}
\usepackage{amsfonts}
\usepackage{amssymb}
\newtheorem{thm}{Theorem}[section]

\newtheorem{cor}[thm]{Corollary}

\newtheorem{remark}[thm]{Remark}

\newtheorem{prop}[thm]{Proposition}
\newtheorem{defn}[thm]{Definition}
\newtheorem{ex}[thm]{Example}
\newcommand{\bb}[1]{\mathbb{#1}}
\newcommand{\cl}[1]{\mathcal{#1}}

\begin{document}
\title[Quasimultipliers of Operator Spaces]
{Quasimultipliers of Operator Spaces}
\author[Masayoshi Kaneda]{Masayoshi Kaneda}
\address{Masayoshi Kaneda: Department of Mathematics, 103 
Multipurpose Science
and Technology Building, University of California, Irvine, 
Irvine, CA
92697-3875 U.S.A.}
\email{mkaneda@math.uci.edu, URL: 
http://www.math.uci.edu/$\sim$mkaneda/}
\author[Vern I. Paulsen]{Vern I. Paulsen}
\address{Vern I. Paulsen: Department of Mathematics, University 
of Houston,
4800 Calhoun Road, Houston, TX 77204-3008 U.S.A.}
\email{vern@math.uh.edu, URL: http://www.math.uh.edu/$\sim$vern/}
\thanks{Research supported in part by a grant from the National 
Science
Foundation}
\begin{abstract}We use the injective envelope to study 
quasimultipliers of
operator spaces. We prove that all representable operator 
algebra products that
an operator space can be endowed with are induced by 
quasimultipliers. We
obtain generalizations of the Banach-Stone theorem.
\end{abstract}
\maketitle
\section{Introduction}
We begin with some general algebraic comments, inspired by 
\cite{BMS}, that
make clear the role that quasimultipliers can play. Let $\cl A$ 
be an
algebra and let $X\subseteq\cl A$ be a subspace. We shall call 
an element $z\in
\cl A$ a {\em quasimultiplier of X (relative to $\cl A$)} 
provided that $XzX 
\subseteq X$, i.e., $x_1zx_2\in X$ for every $x_1, x_2 \in X$. 
Clearly, the set
of quasimultipliers of $X$ is a linear subspace of $\cl A$. 
Moreover, each
quasimultiplier $z$ induces a bilinear map $m_z: X\times X\to X$ 
defined by
$m_z(x_1,x_2):=x_1zx_2$. The associativity of the product on 
$\cl A$, implies
that each $m_z$ is an associative bilinear map and hence can be 
regarded as a
product on $X$. This product gives $X$ the structure of an 
algebra, which we
denote by $(X, m_z)$. There are two homomorphisms, $\pi_l$ and 
$\pi_r$ from
this algebra into $\cl A$, defined by $\pi_l(x):=xz$ and 
$\pi_r(x):=zx$.

The range of $\pi_l$ is contained in the subalgebra of {\em left 
multipliers of
X (relative to $\cl A$)}, while the range of $\pi_r$ is 
contained in the
subalgebra of {\em right multipliers of X (relative to $\cl 
A$)}. Recall that
an element $a\in\cl A$ is called a {\em left (respectively, 
right) multiplier}
of $X$ provided that $aX\subseteq X$ (respectively, $Xa\subseteq 
X$). Finally,
note that since the quasimultipliers are a linear subspace of 
$\cl A$, the set
of ``products" on $X$ that one obtains in this manner is a 
linear subspace of
the vector space of bilinear maps from $X\times X$ into $X$.

In general, linear (or, even convex) combinations of associative 
bilinear maps
need not be associative. For an example of this phenomena, 
consider $X=\bb C^2$
and the associative bilinear maps, $m_1((a,b), (c,d)):=(ac, bd)$ 
and $m_2((a,
b), (c,d)):=(ac, bc)$. Their convex combination, 
$m:=(m_1+m_2)/2$, is not
associative.

One shortcoming of the above representation of quasimultipliers 
is that it is
{\em extrinsic}. The quasimultipliers that one obtains and their 
induced
bilinear maps, could easily depend on the algebra $\cl A$ and on 
the particular
embedding of $X$ into $\cl A$ and not on {\em intrinsic} 
properties of $X$.
Thus, the totality of bilinear maps that one could obtain in 
this manner would
be a union of linear spaces, taken over all embeddings of $X$ 
into an algebra,
which would no longer need to be a linear space.

In this paper, we develop a theory of quasimultipliers of 
operator spaces and
then use the injective envelope to give an intrinsic 
characterization of
quasimultipliers and of their associated bilinear maps. Among 
the results that
we obtain are that an operator space endowed with a completely 
contractive
product can be represented completely isometrically as an 
algebra of operators
on some Hilbert space if and only if the product is a bilinear 
map that belongs
to this space of ``bilinear quasimultipliers". As a corollary, 
we find that the
set of ``representable" completely contractive products is a 
convex set. In
fact, it is affinely isomorphic with the unit ball of the space 
of
quasimultipliers.

We then turn our attention to generalizations of the 
Banach-Stone theorem. Our 
basic result is that a linear complete isometry between any two 
operator
algebras induces a quasimultiplier and that by using the 
quasimultiplier, one
recovers earlier generalizations of the Banach-Stone theorem. 
\section{Quasimultipliers}
In this section we introduce various spaces of quasimultipliers 
of an operator
space and develop some of their key properties. Let $X$ be an 
operator space,
$\cl H$ be a Hilbert space, $B(\cl H)$ denote the algebra of 
bounded, linear
operators on $\cl H$, and let $\phi: X\to B(\cl H)$ be a 
complete isometry. We set $$QM_{\phi}(X):=\{z\in B(\cl H):
\phi(X)z\phi(X)\subseteq\phi(X)\},$$ and we call
$QM_{\phi}(X)$ the space of {\em quasimultipliers of $X$ 
relative to $\phi$}.
Note that $QM_{\phi}(X)$ is a norm closed subspace of $B(\cl H)$.

Each $z\in QM_{\phi}(X)$ induces a bilinear map, $m_z: X\times 
X\to X$ defined
by $m_z(x_1, x_2):=\phi^{-1}(\phi(x_1)z\phi(x_2))$. The bilinear 
map $m_z$ is
completely bounded in the sense of Christensen-Sinclair, that 
is, its linear
extension is completely bounded as a map from $X\otimes_h X$ to 
$X$ with
$\|m_z\|_{cb}\le\|z\|.$ We say that $m_z$ is {\em the bilinear 
map induced by
z}.
\begin{defn} Let $QMB(X)$ denote the set of bilinear maps from 
$X\times X$ to
$X$, that are of the form $m_z$ for some quasimultiplier $z$ and 
some
completely isometric map $\phi$ from $X$ into the bounded linear 
operators on
some Hilbert space. For $m\in QMB(X)$ we set 
$\|m\|_{qm}:=\inf\{\|z\|:
m=m_z\}$, where the infimum is taken over all possible 
completely isometric
maps $\phi$ and $z$ as above.
\end{defn}
Note that by the above remarks, every $m\in QMB(X)$ is 
completely bounded as a
bilinear map and $\|m\|_{cb}\le\|m\|_{qm}$. We shall show in
Example~\ref{ex: sharp} that this inequality can be sharp. For a 
fixed map
$\phi$, the set of bilinear maps, $\{m_z: z\in QM_{\phi}(X)\}$ 
is a linear
subspace of the set $QMB(X)$.  However, since $QMB(X)$ is the 
union of these
subspaces, it is not clear that it is a vector subspace of the 
vector space of
bilinear maps from $X\times X$ to $X$. We shall prove that it is 
a vector space
later.

The above definitions are extrinsic, in the sense that they 
could depend on the
particular embedding. We now seek intrinsic characterizations of 
these maps by
using the injective envelope as in \cite{FP} and \cite{BP} (see 
also
\cite{Pa}).

We begin by recalling a construction used in \cite{BP}, but we 
prefer the
notation from \cite{Pa}.

Recall that if $X\subseteq B(\cl K, \cl H)$ is a (concrete) 
operator space,
then we may form the (concrete) operator system, in $B(\cl 
H\oplus\cl K)$,$$\cl
S_X:=\left\{\left(\begin{array}{c|c}\lambda I_{\cl H}&x\\
\multispan2\hrulefill\\ y^*&\mu I_{\cl 
K}\end{array}\right)\colon \ \lambda,
\mu\in\bb C, x, y\in X\right\}.$$Given a complete isometry 
$\varphi\colon \
X\to B(\cl K_1, \cl H_1)$, the operator system$$\cl
S_{\varphi(X)}:=\left\{\left(\begin{matrix}\lambda I_{\cl 
H_1}&\varphi(x)\\
\varphi(y)^*&\mu I_{\cl K_1}\end{matrix}\right)\colon \ \lambda, 
\mu\in\bb C,
x, y\in X\right\}$$is completely order isomorphic to $\cl S_X$ 
via the map,
$\Phi\colon \ \cl S_X\to\cl S_{\varphi(X)}$ defined by
$$\Phi\left(\left(\begin{matrix}\lambda I_{\cl H} & x\\y^*&\mu 
I_{\cl
K}\end{matrix}\right)\right):=\left(\begin{matrix}\lambda I_{\cl
H_1}&\varphi(x)\\\varphi(y)^*&\mu I_{\cl 
K_1}\end{matrix}\right).$$Thus, the
operator system $\cl S_X$ only depends on the operator space 
structure of $X$
and not on any particular representation of $X$.

Since $\bb C\oplus\bb 
C\cong\left\{\left(\begin{smallmatrix}\lambda I_{\cl
H}&0\\0 & \mu I_{\cl K}\end{smallmatrix}\right)\colon \ \lambda, 
\mu\in\bb
C\right\}$ is a $C^*$-subalgebra of $\cl S_X$, $\bb C\oplus\bb 
C$ will still be
a $C^*$-subalgebra of the $C^*$-algebra, $I(\cl S_X)$ with
$\left(\begin{smallmatrix}I_{\cl H} & 0\\0 & 
0\end{smallmatrix}\right)$ and
$\left(\begin{smallmatrix}0 & 0\\0 & I_{\cl 
K}\end{smallmatrix}\right)$
corresponding to orthogonal projections $e_1$ and $e_2$, 
respectively, in the
$C^*$-algebra $I(\cl S_X)$. We have that $e_1+e_2$ is equal to 
the identity and
$e_1\cdot e_2=0$.

A few words on such a situation are in order. Let $\cl A$ be any 
unital
$C^*$-algebra with orthogonal projections $e_1$ and $e_2$ 
satisfying
$e_1+e_2=1$, $e_1\cdot e_2=0$ and let $\pi\colon \cl A\to B(\cl 
H)$ be a
one-to-one unital $*$-homomorphism. Setting $\cl H_1=\pi(e_1)\cl 
H$, $\cl
H_2=\pi(e_2)\cl H$ we have that $\cl H=\cl H_1\oplus\cl H_2$ and 
relative to
this decomposition every $T\in B(\cl H)$ has the form 
$T=(T_{ij})$ where
$T_{ij}\in B(\cl H_j, \cl H_i)$. In particular, identifying $\cl 
A$ with
$\pi(\cl A)$ we have that $$\cl 
A=\left\{\left(\begin{matrix}a_{11} & a_{12}\\
a_{21} & a_{22}\end{matrix}\right)\colon \ a_{ij}\in\cl 
A_{ij}\right\}$$where
$\cl A_{ij}=e_i\cl Ae_j$, with $\cl A_{ii}\subseteq B(\cl H_i)$ 
unital
$C^*$-subalgebras and $\cl A_{21}=\cl A^*_{12}$. The operator 
space $\cl
A_{12}\subseteq B(\cl H_2, \cl H_1)$ will be referred to as a 
{\em corner of\/}
$\cl A$. Note that $\cl A_{11}\cdot\cl A_{12}\cdot\cl 
A_{22}\subseteq\cl
A_{12}$ so that $\cl A_{12}$ is an $\cl A_{11}-\cl 
A_{22}$-bimodule.

Returning to $I(\cl S_X)$, relative to $e_1$ and $e_2$, we wish 
to identify
each of these 4 subspaces. Note that $X\subseteq e_1 I(\cl 
S_X)e_2=I(\cl
S_X)_{12}$. As shown in Chapter 16 of \cite{Pa}, we may 
identify $I(\cl
S_X)_{12}=I(X)$.

We define, $I_{11}(X):=I(\cl S_X)_{11}$ and $I_{22}(X):=I(\cl 
S_X)_{22}$. Thus
we have the following picture of the $C^*$-algebra $I(\cl S_X)$, 
namely,$$I(\cl
S_X)=\left\{\left(\begin{matrix}a & z\\w^* & 
b\end{matrix}\right)\colon \ a\in
I_{11}(X), b\in I_{22}(X), z, w\in I(X)\right\}$$where 
$I_{11}(X)$ and
$I_{22}(X)$ are injective $C^*$-algebras and $I(X)$ is an 
operator
$I_{11}(X)-I_{22}(X)$-bimodule. Moreover, the fact that $I(\cl 
S_X)$ is a
$C^*$-algebra means that for $z, w\in 
I(X)$,$$\left(\begin{matrix}0 & z\\w^* &
0\end{matrix}\right)\left(\begin{matrix}0 & z\\w^* &
0\end{matrix}\right)=\left(\begin{matrix}zw^* & 0\\0 &
w^*z\end{matrix}\right)$$and consequently there are natural 
products, $z\cdot
w^*\in I_{11}(X)$, $w^*\cdot z\in I_{22}(X)$.

It is interesting to note that setting $\langle z, 
w\rangle=zw^*$ defines an
$I_{11}(X)$-valued inner product that makes $I(X)$ a Hilbert 
$C^*$-module over
$I_{11}(X),$ but we shall not use this additional structure.
\begin{defn}We set $QM(X):=\{z\in I(X)^*: XzX\subseteq X\}$, 
where the products
are all taken in the $C^*$-algebra $I(S_X)$.
\end{defn}
In \cite{BMS} a theory was developed of quasimultipliers of 
Hilbert
$C^*$-bimodules. If $X$ is a Hilbert $C^*$-bimodule, then they 
also use the
notation, $QM(X),$ for their quasimultipliers. We warn the 
reader, and
apologize, that although we are using the same notation, our 
quasimultipliers
and theirs are not the same objects. In the first place their 
quasimultipliers
are defined using the second dual of the linking algebra and 
their
quasimultipliers are a subset of the second dual of $X$.

Briefly, if $X$ is a Hilbert $A-B$-bimodule, then an element $t$ 
in the second
dual of $X$ is a quasimultiplier in the sense of \cite{BMS} 
provided that
$AtB\subseteq X$ where the products are defined in the second 
dual of the
linking algebra. Note that if $A$ and $B$ are both unital 
$C^*$-algebras, then
this forces $t\in X$.

To make our quasimultiplier theory fit with theirs a bit better, 
we should have
perhaps taken adjoints of elements so that our $QM(X)$ is a 
subset of $I(X)$.
However, this alternate definition would have made several 
natural maps, that
we define later, conjugate linear and we believe that it would 
have led to the
``opposite", i.e., transposed operator space structure.
\begin{thm}\label{th: situation}Let $X$ be an operator space, 
$\cl H$ be a Hilbert 
space and let
$\phi: X \to B(\cl H)$, be a complete isometry. Then there 
exists a unique,
completely contractive map $\gamma: QM_{\phi}(X)\to QM(X)$ such 
that
$\phi(x_1)z\phi(x_2)=\phi(x_1\gamma(z) x_2)$ for all $x_1, 
x_2\in X$ and every
$z\in QM_{\phi}(X)$.
\end{thm}
\begin{proof}The proof is similar to that of Theorem~1.7 in 
\cite{BP}. Let $\cl
S_{\phi(X)}\subseteq B(\cl H\oplus\cl H)$ be the concrete 
operator system
defined above and let $C^*(S_{\phi(X)})$ be the $C^*$-subalgebra 
of $B(\cl
H\oplus \cl H)$ that it generates. The $C^*$-subalgebra of 
$I(S_X)$ generated
by $S_X$ is known to be the $C^*$-envelope of $S_X, C^*_e(S_X)$. 
Consequently,
by \cite{H2} Corollary~4.2 the identity map on $S_X$ extends to 
be a
surjective $*$-homomorphism $\pi: C^*(S_{\phi(X)})\to 
C^*_e(S_X)$.

Let $\Gamma: B(\cl H\oplus\cl H)\to I(S_X)$ be a completely 
positive map that
extends this $*$-homomorphism. Since $\Gamma$ extends $\pi$, it 
will be a
$\pi$-bimodule map, that is, for $A, B\in C^*(S_{\phi(X)})$ we 
will have that
$\Gamma(ATB)=\pi(A)\Gamma(T)\pi(B)$. This forces $\Gamma$ to be 
a matrix of
maps, that is, for $T=\begin{pmatrix}T_{11} & T_{12}\\T_{21} & 
T_{22}
\end{pmatrix}\in B(\cl H\oplus\cl H)$ we will have that
$\Gamma(T)=\begin{pmatrix}\gamma_{11}(T_{11}) & 
\gamma_{12}(T_{12})\\
\gamma_{21}(T_{21}) & \gamma_{22}(T_{22})\end{pmatrix}$.

In particular, for $x\in X$ and $z\in QM_{\phi}(X)$ we will have 
that
$\Gamma(\begin{pmatrix}0 & \phi(x)\\z & 
0\end{pmatrix})=\begin{pmatrix}0 & x\\
\gamma_{21}(z) & 0\end{pmatrix}$.

Hence by the $\pi$-bimodule property, for $x_1, x_2\in X$ and 
$z\in
QM_{\phi}(X)$ we have \begin{multline*}\Gamma(\begin{pmatrix}0 &
\phi(x_1)z\phi(x_2)\\0 & 0\end{pmatrix})=\Gamma(\begin{pmatrix}0 
& \phi(x_1)\\0
& 0\end{pmatrix}\begin{pmatrix}0 & 0\\z & 
0\end{pmatrix}\begin{pmatrix}0 &
\phi(x_2)\\0 & 0\end{pmatrix})\\=\begin{pmatrix}0 & x_1\\0 & 
0\end{pmatrix}
\begin{pmatrix}0 & 0\\\gamma_{21}(z) & 
0\end{pmatrix}\begin{pmatrix}0 & x_2\\0
& 0\end{pmatrix}=\begin{pmatrix}0 & x_1\gamma_{21}(z)x_2\\0 & 
0\end{pmatrix}.
\end{multline*}

Thus, it follows that $\gamma_{21}(z)\in QM(X)$ and 
$\phi(x_1\gamma_{21}(z)
x_2)=\phi(x_1)z\phi(x_2)$. Because $\Gamma$ is a unital, 
completely positive
map, $\gamma_{21}$ is completely contractive.

Finally, the uniqueness of the map comes from the following 
observation.
Suppose $q_1, q_2\in QM(X)$ have the property that 
$x_1q_1x_2=x_1q_2x_2$ for
every $x_1, x_2 \in X$. This implies that $(x_1q_1-x_1q_2)X=0$ 
and so by
Corollary~1.3 of \cite{BP}, we have that $x_1(q_1-q_2)=0$ for 
every $x_1\in X$.
Applying the corollary again, we see that 
$X(q_1-q_2)(q_1-q_2)^*=0$ and so
$q_1=q_2$.
\end{proof}
\begin{remark}Let $A$ be a $C^*$-algebra and let $\pi_{u}: A\to 
B(\cl H_u)$ be
its universal representation. The classical definition of the 
quasimultiplier
space of $A$ as given in \cite{Pe} is the set $QM_{\pi_u}(A)$. 
In \cite{Ka},
the first author proves that the map $\gamma: QM_{\pi_u}(A)\to 
QM(A)$ is an
onto complete isometry, and preserve quasimultiplication. Thus, 
at least in the
case of a $C^*$-algebra our definition and the classical 
definition agree.
\end{remark}
\begin{cor}Let $X$ be an operator space. The map $z\mapsto m_z$ 
from $QM(X)$ to
$QMB(X)$ equipped with $\|\cdot\|_{qm}$ is an onto isometry, 
where $m_z:
X\times X\to X$ is defined by $m_z(x_1, x_2):=x_1zx_2$. 
Consequently,
$QMB(X)$ is a linear subspace of the vector space of bilinear 
maps from
$X\times X$ to $X$. 
\end{cor}
Let $X$ be an operator space. Given an associative, bilinear map 
$m: X\times
X\to X$, we let $(X, m)$ denote the resulting algebra. We let 
$CCP(X)$ denote
the set of associative bilinear maps on $X$ that are completely 
contractive in
the sense of Christensen-Sinclair, that is, $CCP(X)$ denotes the 
set of {\em
completely contractive products on X}. We let $OAP(X)\subseteq 
CCP(X)$ denote
those maps such that the algebra $(X, m)$ has a completely 
isometric
homomorphism into $B(\cl H)$ for some Hilbert space $\cl H$. 
That is, $OAP(X)$
denotes the set of {\em operator algebra products on X}.

We let $CBP(X)$ denote those associative, bilinear maps from 
$X\times X$ to $X$
(i.e., products on $X$), that are completely bounded in the 
sense of
Christensen-Sinclair, that is the set of {\em completely bounded 
algebra
products}. By a result of \cite{Bl}, $m\in CBP(X)$ if and only 
if there exists
a Hilbert space $\cl H$ and a completely bounded homomorphism, 
$\pi: (X, m)\to
B(\cl H)$ with completely bounded inverse, $\pi^{-1}: 
\pi(X)\to(X, m)$. In
fact, by \cite{Bl} one may choose $\pi$ such that
$\|\pi\|_{cb}\|\pi^{-1}\|_{cb}\le 1+\epsilon$ for any 
$\epsilon\ge0$. Finally,
we let $SOAP(X)\subseteq CBP(X)$ denote those associative, 
bilinear maps for
which one can choose, $\pi$ satisfying, 
$\|\pi\|_{cb}\|\pi^{-1}\|_{cb}=1$,
i.e., such that $\pi$ is a scalar multiple of a complete 
isometry. These are
the {\em scaled operator algebra products}.

The following theorem illustrates the importance of 
quasimultipliers.
\begin{thm}\label{th: oap}Let $X$ be an operator space. Then 
$QMB(X)=SOAP(X)$
and $OAP(X)=\{m\in QMB(X): \|m\|_{qm}\le1\}$.
\end{thm}
\begin{proof}We prove the second equality first. Let $m\in 
OAP(X)$ and let
$\pi: (X, m)\to B(\cl H)$ be a completely isometric 
homomorphism. Then $I_{\cl
H}$ is a quasimultiplier of $\pi(X)$ that induces the bilinear 
map $m$. Thus,
$m\in QMB(X)$ and $\|m\|_{qm}\le1$. Conversely, let $m\in 
QMB(X)$ with
$\|m\|_{qm}\le1$. Then there exists $z\in QM(X)$ with 
$\|z\|\le1$ such that
$m=m_z$.

As \cite{BRS} Remark 2, define $\pi: (X, m)\to I(S_X)$ by
$$\pi(x):=\begin{pmatrix}xz & x\sqrt{1-zz^*}\\
0 & 0\end{pmatrix}.$$It is easily seen that,
$\pi(x_1)\pi(x_2)=\pi(x_1zx_2)=\pi(m(x_1,x_2))$ and that
$\|\pi(x)\|^2=\|\pi(x)\pi(x)^*\|=\|xx^*\|=\|x\|^2$, where all 
products take
place in $I(S_X)$. Thus, $\pi$ is an isometric homomorphism. The 
proof that
$\pi$ is completely isometric is similar and thus $m\in OAP(X)$.

If $m\in SOAP(X)$ and $\pi: (X,m)\to B(\cl H)$ is a completely 
bounded
homomorphism such that $\phi=\pi/\|\pi\|_{cb}$, is a complete 
isometry, then
$I_{\cl H}$ is a quasimultiplier of $\phi(X)$. We have 
that$$\phi(x_1)I_{\cl
H}\phi(x_2)=\|\pi\|_{cb}^{-2}\pi(x_1)\pi(x_2)=\|\pi\|_{cb}^{-2}\pi(m(x_1,
x_2))=\|\pi\|_{cb}^{-1}\phi(m(x_1, x_2)).$$Hence, $m\in QMB(X)$ 
and
$\|m\|_{qm}\le\|\pi\|_{cb}$.

Conversely, if $m(x_1, x_2)=x_1zx_2$ for $z\in QM(X)$ with 
$\|z\|=r$, then
consider$$\pi(x)=\begin{pmatrix}xz & x\sqrt{r^2-zz^*}\\0 & 
0\end{pmatrix}$$and
argue as above to prove that $\pi$ is a homomorphism and is $r$ 
times a
complete isometry.
\end{proof}
\begin{cor}Let $X$ be an operator space, then $OAP(X)$ is a 
convex set, and
$SOAP(X)$ is a vector space.
\end{cor}
The following example shows that, in general, $CCP(X)$ and 
$CBP(X)$ are not
convex sets.
\begin{ex}\label{ex: nonconvex}The following example shows that, 
in general,
$OAP(X)\ne CCP(X)$, $SOAP(X)\ne CBP(X)$, and that the sets 
$CCP(X)$ and
$CBP(X)$ need not be convex.

Let $X=C_2$, that is, the subspace of the two by two matrices, 
$M_2$ consisting
of those matrices that are 0 in the second column and whose 
first column is
arbitrary. Let $m_1(\begin{pmatrix}a\\b\end{pmatrix}, 
\begin{pmatrix}c\\d
\end{pmatrix}):=\begin{pmatrix}ac\\bd\end{pmatrix}$ and 
$m_2(\begin{pmatrix}a\\
b\end{pmatrix}, 
\begin{pmatrix}c\\d\end{pmatrix}):=\begin{pmatrix}ac\\bc
\end{pmatrix}$ be as in the introduction.

Since $m_2$ is the product on $C_2$ induced by the inclusion of 
$C_2$ into
$M_2$, we have that $m_2\in OAP(C_2)\subseteq CCP(C_2)$.

We claim that $m_1\in CCP(C_2)$. To see this claim, note that if 
we identify
$M_n(C_2)=\{\begin{pmatrix}A\\B\end{pmatrix}: A, B\in M_n\}$, 
then
$m_1^{(n)}(\begin{pmatrix}A\\B\end{pmatrix}, \begin{pmatrix}C\\D
\end{pmatrix})=\begin{pmatrix}AC\\BD\end{pmatrix}$. If 
$\|\begin{pmatrix}A\\B
\end{pmatrix}\|\le1$ and 
$\|\begin{pmatrix}C\\D\end{pmatrix}\|\le1$, then
$\|\begin{pmatrix}AC\\BD\end{pmatrix}\|^2=\|C^*A^*AC+D^*B^*BD\|$. 
However,
$0\le C^*A^*AC+D^*B^*BD\le C^*C+D^*D\le I$ and so 
$\|m_1\|_{cb}\le1$. Thus,
$m_1\in CCP(X)$ as claimed.

Since $(m_1+m_2)/2$ is not even associative, we see that 
$CCP(C_2)$ is not
convex and hence cannot be equal to $OAP(C_2)$. Since $m_1$ and 
$m_2$ are both also in $CBP(C_2)$, we have that this set is 
also not convex and hence it is not equal to $SOAP(C_2).$
Also, since $SOAP(C_2)$ is convex and $m_2 \in SOAP(C_2)$, it 
must be the case that $m_1$ is not in $SOAP(C_2).$

This last fact can also be seen by using the injective envelope 
and
Theorem~\ref{th: oap}. Since $C_2$ is already an injective 
operator space, we
have that $I(C_2)^*=C_2^*=R_2$, where $R_2$ denotes the 
corresponding row
space. Since $C_2R_2C_2\subseteq C_2$, we have that 
$QM(C_2)=R_2$. Now it is
easily checked that there is no $z=(e, f)\in R_2$ such that 
$m_1(x_1,
x_2)=x_1zx_2$ and so $m_1$ is not in $QMB(C_2)=SOAP(C_2)$.

Finally, note that $\pi: (C_2, m_1)\to M_2$ defined by 
$\pi(\begin{pmatrix}a\\b
\end{pmatrix}):=\begin{pmatrix}a & 0\\0 & b\end{pmatrix}$ is a 
completely
contractive homomorphism with completely bounded inverse. Which 
gives a direct way, independent of Blecher's theorem \cite{Bl}, 
to see that $(C_2,m_1)$ is completely boundedly representable.

It is also interesting to note that for the natural inclusion 
$\phi: C_2\to
M_2$, we have that $QM_{\phi}(C_2)$ is all of $M_2$, while 
$QM(C_2)= R_2$ is
2-dimensional. Thus, we see that the map $\gamma$ of 
Theorem~\ref{th: situation} need not be one-to-one.
\end{ex}
\begin{remark}Although, $OAP(X)$ is a convex set and $CCP(X)$ is 
not in
general, we know little else about the structure of $CCP(X)$ or 
about the
subset of $CCP(X)$ consisting of those products that can be 
induced by a
completely contractive, but not completely isometric 
representation. The
product $m_1$ from Example~\ref{ex: nonconvex} is one such 
product.

For a finite dimensional vector space the set of associative
bilinear maps is an
algebraic set. Thus, when $X$ is a finite-dimensional 
operator space,
$CCP(X)$ is the intersection of this algebraic set with the 
set of
completely contractive bilinear maps. Generally, the set of 
completely contractive bilinear maps need not even be a 
semialgebraic set. But it is still possible that the 
intersection, $CCP(X)$, is a semialgebraic set.
\end{remark}
The following result illustrates how some of the results of 
\cite{BK} on
operator algebras with one-sided identities can be deduced from 
the theory of
quasimultipliers.
\begin{prop}\label{pr: bk}Let $X$ be an operator space and let 
$m\in OAP(X)$.
Then $(X, m)$ has a right contractive identity $e$ if and only 
if $m=m_z$ where
$z\in QM(X)$ satisfies $z^*=e$, and $zz^*$ is the identity of 
$I_{22}(X)$. In
this case the map $x\to xz$ defines a completely isometric 
homomorphism of $(X,
m)$ into $M_l(X)$.
\end{prop}
\begin{proof}Since $m\in OAP(X)$ we have that $m=m_z$ for some 
$z\in QM(X)$
with $\|z\|\le1$.

Assume that $(X, m)$ has a contractive right identity $e$. Then 
we have that
for every $x\in X, x=m(x, e)=xze$. Hence, $x(1_{22}-ze)=0$ for 
every $x\in X$.
By \cite{BP} Corollary 1.3, this implies that $1_{22}-ze=0$. But 
since both 
$z$ and $e$ are
contractions, $e=z^*$ must hold.

Conversely, if $z^*=e$ and $zz^*=1_{22}$, then clearly, $m(x, 
e)=xze=x$ and so
$e$ is a contractive, right identity.

Finally, since $\|x\|=\|xzz^*\|\le\|xz\|$, we see that the 
completely
contractive homomorphism of $(X, m)$ into $M_l(X)$ given by 
$x\to xz$ is a
complete isometry.
\end{proof}

Note that by the above result, the relationship between the 
product $m$ and the
product in $I(S_X)$ is that $m(x_1, x_2)=x_1e^*x_2$.

There is an analogous result for left identities.

It is possible for a concrete algebra of operators to have a 
two-sided identity
$e$ of norm greater than one. For an example see \cite{Pa}, 
page 279. 
In this case the
multiplication will still be given by a contractive 
quasimultiplier $z$ and one
has $ze=1_{22}, ez=1_{11}$ but one no longer has that $e=z^*$.

We close this section with a number of examples of spaces of 
quasimultipliers
that illustrate the limits of some of the above results.
\begin{ex}\label{ex: sharp}This example shows that it is 
possible to have
$\|z\|>\|m_z\|_{cb}$ for a quasimultiplier. It is based on 
Example~4.4 of
\cite{B2}.

Let $\cl A\subseteq M_3$ denote the subalgebra that is the span 
of $\{E_{12},
E_{13}, I_3\},$ where $I_3$ denotes the identity matrix. Let 
$Q:=I_3+J$ where
$J$ is the matrix whose entries are all 1's, and set 
$P:=Q^{1/2}$. A little
calculation shows that $P=I_3+\frac{1}{3}J$ and that 
$P^{-1}=I_3-\frac{1}{6}J$.

We let $\cl X=\cl AP$. Since the $C^*$-subalgebra of $M_3$ 
generated by $\cl X$
is all of $M_3$, which is irreducible, one finds that 
$I_{11}(\cl X)=I(\cl
X)=M_3,$ with the usual product. From this we see quite easily  
that $M_l(\cl
X)=\cl A$ and $QM(\cl X)=P^{-1}\cl A$.

Let $Z:=P^{-1}(E_{12}- E_{13})$, so that $\|Z\|=\sqrt{3/2}$. 
Writing 
$X_i, Y_i\in\cl
X$ as $X_i=(a_iI_3+N_i)P, Y_i=(b_iI_3+M_i)P$ where $N_i$ and 
$M_i$ are in the
span of $E_{12}$ and $E_{13}$, we have that, $m_Z(X_i,
Y_i):=X_iZY_i=a_ib_i(E_{12}-E_{13})P$.

Since $\|(E_{12}-E_{13})P\|=\sqrt{2}$, we have that 
$\|m_Z\|=\kappa\sqrt{2}$, where
$\kappa=\sup\{|\sum a_ib_i|: \|(X_1, \ldots, X_n)\|\le1, 
\|(Y_1^*, \dots,
Y_n^*)\|\le1\}$ and $m_Z$ is the linear map $m_Z: \cl A 
\otimes_h \cl A \to \cl A.$ The first inequality implies that 
$I_3\ge\sum 
X_iX_i^*$.
Examining the (3,3)-entry of these matrices leads to the 
conclusion that
$2\sum|a_i|^2\le1$. The second inequality implies that
$\sum (b_iI_3 + M_i)^*(b_iI_3 + M_i) \le P^{-2}.$ Examining the 
 (1,1)-entry of these matrices leads to the conclusion that 
$\sum|b_i|^2\le 3/4$ and hence, $\kappa\le \sqrt{3/8}$.

Thus, we are led to conclude that 
$\|m_Z\|\le\sqrt{3/4}<\sqrt{3/2}$. 

The same calculation, using 
matrix coefficients for the entries of $X_i$ and $Y_i$ shows 
that $\|m_Z\|_{cb}\le\sqrt{3/4}$ too, 
and so the
result follows.  
Indeed, if we write $X,Y \in M_n(\cl A)$ as $X= A \otimes 
E_{11} + C \otimes E_{12} + D \otimes E_{13} + A \otimes E_{22} 
+ A \otimes E_{33}$ and $Y= B \otimes E_{11} + G \otimes E_{12} 
+ H \otimes E_{13} + B \otimes E_{22} + B \otimes E_{33}$ where 
$A,B,C,D,G,H$ are in $M_n$, then it is easily seen that
$m_Z^{(n)}(X,Y) = AB \otimes (E_{12} -E_{13})P$
Now taking $X_i \in M_n(A)$ with $\|(X_1, \ldots, X_m)\| \le
1$ implies that $\sum X_iX_i^* \le I_3 \otimes I_n$  and 
examining in
the
(1,1)-block of this matrix inequality yields,
$\sum A_iA_i^* \le I_n$.
Similarly one gets that $\sum B_i^*B_i \le 3/4I_n$
Hence, $\sum A_iB_i = (A_1, \ldots, A_m) \cdot (B_1, \ldots, 
B_m)^t$
and these estimates yield that the norm of the row of A's is
less than one
and norm of the column of B's is less than $\sqrt{3/4}.$

In fact, it is not too hard to show that
$\|m_Z\|=\|m_Z\|_{cb}=\sqrt{3/4}$.
\end{ex}
\begin{ex}Let $\{E_{ij}\}$ denote the canonical matrix units and 
let
$X=span\{E_{11}, E_{12}, E_{21}, E_{32}\}\subseteq M_{3, 2}$. We
compute $QM(X)$ for this space and illustrate some of its 
properties.

It is not difficult to show that $I(S_X)=\begin{pmatrix}M_3 & 
M_{3, 2}\\M_{2,
3} & M_2\end{pmatrix}=M_5$, with the obvious identifications. To 
see this one
first shows that since $X$ is a $\cl D_3-\cl D_2$-bimodule, 
where $\cl D_n$
denotes the $n\times n$ diagonal matrices, then any completely 
contractive map
$\Phi$ from $M_{3, 2}$ into itself that fixes $X$ must be a 
bimodule map. From
this it follows that $\Phi$ must be given as a Shur product map, 
but then the
fact that $\Phi$ is completely contractive forces $\Phi$ to be 
the identity
map.

Now a direct calculation shows that $QM(X)=span\{E_{12}, 
E_{23}\}\subseteq
M_{2, 3}$, that $M_l(X)=span\{E_{11}, E_{12}, E_{13}, E_{22}, 
E_{33}\}\subseteq
M_3$ and that \hfill \break \hfill $M_r(X)=span\{E_{11}, 
E_{22}\}\subseteq M_2$.

Note that the span of the products $X\cdot QM(X)$ is not dense 
in $M_l(X)$ but
that the span of the products $QM(X)\cdot X$ is all of $M_r(X)$.

For the contractive quasimultiplier $z=E_{12}+E_{23}$, we see 
that the
induced homomorphism $\pi_l(x)=xz$ into $M_l(X)$ is a complete 
isometry, but
that $\pi_r(x)=zx$ is not even one-to-one. For the 
quasimultipliers $E_{12}$
and $E_{23}$ neither $\pi_l$ nor $\pi_r$ is one-to-one.
\end{ex}
\begin{ex}Let $X=span\{E_{11}+E_{32}, E_{21}+E_{33}\}\subseteq 
M_3$.
This space can be identified as a concrete representation of the 
maximum of
$C_2$ and $R_2$, that is, as the least operator space structure 
on $\bb C^2$
that is greater than both $C_2$ and $R_2$.  We will show that 
$QM(X)=(0)$ and
consequently, there are no non-trivial operator algebra products 
on this
operator space, i.e., $OAP(X)=(0)$. However, since the natural 
maps from $X$ to
the concrete operator algebras $C_2\subseteq M_2$ and 
$R_2\subseteq M_2$ are
both complete contractions, we see that there are at least 4 
different products
(up to scaling) in $CCP(X)$ that have completely contractive 
representations
whose inverses are completely bounded.

To see these claims, one first shows that if one regards 
$S_X\subseteq M_6,$
then $C^*(S_X)=I(S_X)$. From this it follows that 
$I_{11}(X)=M_2\oplus\bb C$, $I_{22}(X)=\bb C\oplus M_2$, 
$I(X)=span\{E_{11},
E_{32}, E_{21}, E_{33}\}$
and that $M_l(X)$ and $M_r(X)$ are both just the scalar 
multiples of the
identity. Once these things are seen, it is straightforward to 
check that
$QM(X)=(0)$.

To prove that $C^*(S_X)= I(S_X)$, first note that there is a 
*-homomorphism of $C^*(S_X)$ onto the
C*-subalgebra of $I(S_X)$
generated by the copy of $S_X,$ i.e., onto the boundary
C*-algebra. But the
original C*-algebra has only 2 ideals that could be the
kernel of this
map. Now argue that if you mod out by either ideal then you
will not have
a 2-isometry on $S_X.$ Hence this homomorphism must be 1-1. But 
$C^*(S_X)$ is
injective so we are done.

\end{ex}
\section{A Non-commutative Banach-Stone Theorem}In this section 
we use
quasimultipliers to obtain a characterization of linear complete 
isometries
from one operator algebra onto another. Our theorem needs no 
assumptions
concerning the existence of units or approximate units.

To understand the statement of the theorem, it is perhaps 
instructive to keep
the following example in mind. Let $\cl A\subseteq B(\cl H)$ be 
a unital
$C^*$-subalgebra and let $\cl B:=\{\begin{pmatrix}0 & a\\0 & 
0\end{pmatrix}:
a\in\cl A\}\subseteq B(\cl H\oplus\cl H)$. These are both 
algebras of
operators, although the product of any two elements of $\cl B$ 
is 0. The
identification of $\cl A$ with $\cl B$ is a complete isometry, 
onto, but
clearly the only possible homomorphism between these two 
algebras is the 0 map.
However, in this example one sees that the left multipliers of 
$\cl B$ can be
identified with the $C^*$-algebra $\cl A$. 
\begin{thm}Let $\cl A$ and $\cl B$ be algebras of operators and 
let $\psi: \cl
A\to\cl B$ be a linear complete isometry that is onto. Then we 
have the
following:
\begin{enumerate}
\item there exists a unique $z\in QM(\cl A)$, with $\|z\|\le1$ 
such that
$\psi(a_1)\psi(a_2)=\psi(a_1za_2)$ for every $a_1, a_2\in\cl A$;
\item there exists a unique $w\in QM(\cl B)$, with $\|w\|\le1$ 
such that
$\psi(a_1a_2)=\psi(a_1)w\psi(a_2)$ for every $a_1, a_2\in\cl A$;
\item setting $\pi_l(a)=\psi(a)w$ and $\pi_r(a)=w\psi(a)$ 
defines completely
contractive homomorphisms of $\cl A$ into $M_l(\cl B)$ and 
$M_r(\cl B)$,
respectively;
\item if $\cl A$ has a contractive right (respectively, left) 
approximate
identity, then $\pi_l$ (respectively, $\pi_r$) is a completely 
isometric
homomorphism;
\item if $\cl A$ has a contractive right (respectively, left) 
identity, $e$,
then $w=\psi(e)^*$ and $ww^*$ is the identity of $I_{22}(\cl B)$ 
(respectively,
$w^*w$ is the identity of $I_{11}(\cl B)$).
\end{enumerate}
\end{thm}
\begin{proof}Set $\gamma:=\psi^{-1}$ and define $m: \cl 
A\times\cl A\to\cl A$
by $m(a_1, a_2):=\gamma(\psi(a_1)\psi(a_2))$. It is easily 
checked that$$m(a_1,
m(a_2, a_3))=\gamma(\psi(a_1)\psi(a_2)\psi(a_3))=m(m(a_1, a_2), 
a_3),$$so that
$m$ is an associative bilinear map on $\cl A$ and defines a new 
product on $\cl
A$. Moreover, because the product on $\cl B$ is completely 
contractive this new
product on $\cl A$ is completely contractive and the map $\psi: 
(\cl A,
m)\to\cl B$ is a completely isometric algebra isomorphism.

Thus, since $\cl B$ is an algebra of operators, we see that 
$m\in OAP(\cl A)$
and hence by Theorem~\ref{th: oap} there exists a unique $z \in 
QM(\cl A)$ such
that $m(a_1, a_2)=a_1za_2$. Hence, $\psi(a_1za_2)=\psi(m(a_1,
a_2))=\psi(a_1)\psi(a_2)$ and (1) 
follows.

Applying (1) to $\gamma$ yields $w\in QM(\cl B)$ such that
$$\gamma(\psi(a_1)w\psi(a_2))=\gamma(\psi(a_1))\gamma(\psi(a_2))=a_1a_2.$$Thus,
$\psi(a_1)w\psi(a_2)=\psi(a_1a_2)$ and so (2) follows.

To see (3), note that
$\pi_l(a_1)\pi_l(a_2)=\psi(a_1)w\psi(a_2)w=\psi(a_1a_2)w=\pi_l(a_1a_2)$ 
with a
similar calculation for $\pi_r$. Since
$\pi_l(a)b=\psi(a)w\psi(\gamma(b))=\psi(a\gamma(b))\in\cl B$ for 
every $b\in\cl
B$, we see that $\pi_l(a)\in M_l(\cl B)$ for every $a\in\cl A,$ 
with a similar
calculation for $\pi_r$.

Now let $\{e_{\alpha}\}$ be a contractive, approximate right 
identity for $\cl
A$. We then have that $\|\psi(a)\|=\lim\|\psi(ae_{\alpha})\|
=\lim\|\psi(a)w\psi(e_{\alpha})\|=\lim\|\pi_l(a)\psi(e_{\alpha})\|
\le\|\pi_l(a)\|$. Thus, $\pi_l$ is an isometry. The proof that 
$\pi_l$ is a
complete isometry and the case for $\pi_r$ are similar.

Finally, if $\cl A$ has a right identity $e$, then 
$\psi(a)=\psi(ae)
=\psi(a)w\psi(e)$. This shows that $bw\psi(e)=b$ for every 
$b\in\cl B$ and
hence $w\psi(e)$ is a right identity for $M_r(\cl B)$. By 
\cite{BP} Corollary
1.3, we have that $w\psi(e)$ is the identity of $I_{22}(\cl B)$. 
Since
$\|w\|\le1$ and $\|\psi(e)\|\le1$, we have that $\psi(e)=w^*$ 
and (5) follows.
\end{proof}
\begin{remark}
\begin{enumerate}
\item When $\cl B$ has a contractive right identity, then one 
may identify $\cl
B\subseteq M_l(\cl B)$, but it is not clear if the image of 
$\pi_l$ maps onto
this copy of $\cl B$. However, in this case it is clear how to 
define a
homomorphism into $\cl B$. Let $\psi(a_0)=e_{\cl B}$ and define 
$\rho: \cl
A\to\cl B$ by setting $\rho(a)=e_{\cl B}w\psi(a)=\psi(a_0a)$. 
Letting the
product in $\cl B$ be denoted by $\odot$ to avoid confusion, we 
have that
$b_1\odot b_2=b_1e_{\cl B}^*b_2$, where the latter product is 
taken in
$I(S_{\cl B})$. Since $e_{\cl B}^*e_{\cl B}=1_{22}$, by 
Proposition~\ref{pr:
bk}, we have that $\rho(a_1)\odot\rho(a_2)
=e_{\cl B}w\psi(a_1)e_{\cl B}^*e_{\cl B}w\psi(a_2)=e_{\cl 
B}w\psi(a_1a_2)
=\rho(a_1a_2)$, and so $\rho$ is a homomorphism. Note that 
$\rho$ is onto $\cl B$
if and only if $a_0\cl A=\cl A$.
\item If one considers $\cl A=\cl B=C_2\subseteq M_2$ and lets 
$\psi$ be the
identity map, then we are in the situation of the last remark. 
Thus, $\pi_l$ is
a complete isometry, but since $M_r(C_2)=\bb C$, we have that 
$\pi_r$ is not a
complete isometry. In fact, it is the compression to the (1, 
1)-entry.
\end{enumerate}
\end{remark}

\section{Further Results on QMB(X)}In \cite{B2} and \cite{BEZ} 
various
characterizations are given of the linear maps of an operator 
space $X$ into
itself that are given as left multiplication by an element from 
the left
multiplier algebra of $X, M_l(X)$. In this section we present 
characterizations
of the bilinear maps of an operator space into itself that are 
in $QMB(X)$.
Among the results that we obtain is a characterization of when 
a linear map
from $X$ into $M_l(X)$ is given as right multiplication by a 
quasimultiplier.
We also identify a subspace of $QM(X)$, related to ternary 
structures on $X$,
that we denote by $TER(X)^*$ for which we have
$\|z\|=\|m_z\|_{qm}=\|m_z\|_{cb}$.

In \cite{BEZ} it was shown that one could determine whether or 
not a linear map
from $X$ into $X$ was given by a contractive left multiplier by 
determining
whether or not an associated linear map was completely 
contractive. The
following is an analogous result for determining when a map is 
given as
multiplication by a quasimultiplier.

Recall that given any operator space $X, R_2(X)$ denotes the 
operator subspace
of $M_2(X)$ consisting of $1\times2$ matrices.

\begin{thm}Let $X$ be an operator space and let $\gamma: X\to 
I_{11}(X)$ be a
linear map. There exists $y\in I(X)^*$ with $\|y\|\le1$ such that
$\gamma(x)=xy$ for every $x\in X$ if and only if the map $\beta: 
R_2(X)\to
I(S_X)$ defined by $\beta((x_1, 
x_2)):=\begin{pmatrix}\gamma(x_1) & x_2\\0 & 0
\end{pmatrix}$ is completely contractive.
\end{thm}
\begin{proof}Note that if such an element $y$ exists, then 
$\beta$ is given, at
least formally, as right multiplication by the matrix 
$\begin{pmatrix}y & 0\\0
& 1_{22}\end{pmatrix}$ and since this matrix has norm 1, $\beta$ 
should be a
complete contraction. To complete this argument, we create a 
$C^*$-algebra
where these products occur.

To this end consider the following $C^*$-algebra,$$\cl 
B:=\begin{pmatrix}
I_{11}(X) & I(X) & I(X)\\I(X)^* & I_{22}(X) & I_{22}(X)\\I(X)^* 
& I_{22}(X) &
I_{22}(X)\end{pmatrix},$$where the products are all induced from 
the products
in $I(S_X)$. Identifying $R_2(X)$ with the subspace 
$\begin{pmatrix}0 & X & X\\
0 & 0 & 0\\0 & 0 & 0\end{pmatrix}\subseteq\cl B$, we see that 
$\beta$ is given
as right multiplication in the $C^*$-algebra $\cl B$ by the 
matrix
$\begin{pmatrix}0 & 0 & 0\\y & 0 & 0\\0 & 1_{22} & 
0\end{pmatrix}$.

For the converse, we must assume that $\beta$ is a complete 
contraction and
produce the element $y$. To this end we create a second 
$C^*$-algebra, $\cl C$
and an operator system $\cl S$.

Let$$\cl C:=\begin{pmatrix}I_{11}(X) & I_{11}(X) & I(X)\\ 
I_{11}(X) & I_{11}(X)
& I(X)\\I(X)^* & I(X)^* & I_{22}(X)\end{pmatrix},$$where the 
products are all
induced from the products in $I(S_X)$ and let$$\cl 
S=\begin{pmatrix}\bb C1_{11}
& X & X\\X^* & \bb C1_{22} & 0\\X^* & 0 & \bb 
C1_{22}\end{pmatrix}\subseteq\cl
B.$$

We define $\Phi: \cl S\to\cl C$ by 
$$\Phi(\begin{pmatrix}\lambda1_{11} & x_1 &
x_2\\x_3^* & \mu1_{22} & 0\\x_4^* & 0 & \nu 
1_{22}\end{pmatrix}):=
\begin{pmatrix}\lambda1_{11} & \gamma(x_1) & x_2\\\gamma(x_3)^* 
& \mu 1_{11} &
0\\x_4^* & 0 & \nu 1_{22}\end{pmatrix}.$$

Since $\beta$ is completely contractive, $\Phi$ is completely 
positive and
since $\cl C$ is clearly an injective $C^*$-algebra, we may 
extend $\Phi$ to a
completely positive map on all of $\cl B$, which we still denote 
by $\Phi$.
Because $\Phi$ fixes the diagonal, it will be a bimodule map 
over the diagonal.
Also note that the compression of $\cl S$ to the span of the 
first and third
entries is a copy of $S_X$ and that $\Phi$ fixes this operator 
system.

By the rigidity properties of the injective envelope, we see 
that necessarily
$$\Phi(\begin{pmatrix}a & 0 & b\\0 & \mu 1_{22} & 0\\c^* & 0 & 
d\end{pmatrix}
=\begin{pmatrix}a & 0 & b\\0 & \mu 1_{11} & 0\\c^* & 0 & 
d\end{pmatrix},$$for
every $a\in I_{11}(X), b, c\in I(X), d\in I_{22}(X)$ and 
$\mu\in\bb C$.

These matrices that are fixed by $\Phi$ form a common 
$C^*$-subalgebra of $\cl
B$ and $\cl C$ and hence $\Phi$ will necessarily be a bimodule 
map over this
$C^*$-subalgebra.

Thus, we will have for any $x\in X$ that 
\begin{multline*}\begin{pmatrix}0 &
\gamma(x) & 0\\0 & 0 & 0\\0 & 0 & 
0\end{pmatrix}=\Phi(\begin{pmatrix}0 & x & 0
\\0 & 0 & 0\\0 & 0 & 0\end{pmatrix})=\Phi(\begin{pmatrix}0 & 0 & 
x\\0 & 0 & 0\\
0 & 0 & 0\end{pmatrix}\cdot\begin{pmatrix}0 & 0 & 0\\0 & 0 & 
0\\0 & 1_{22} & 0
\end{pmatrix})\\=\begin{pmatrix}0 & 0 & x\\0 & 0 & 0\\0 & 0 & 
0\end{pmatrix}
\Phi(\begin{pmatrix}0 & 0 & 0\\0 & 0 & 0\\0 & 1_{22} & 
0\end{pmatrix})
=\begin{pmatrix}0 & 0 & x\\0 & 0 & 0\\0 & 0 & 
0\end{pmatrix}\begin{pmatrix}0 &
0 & 0\\0 & 0 & 0\\0 & y & 0\end{pmatrix},\end{multline*}where 
$y$ is the image
of $1_{22}$ under the restriction of $\Phi$ to the subspaces 
corresponding to
the (3,2)-entries.

Equating entries of the matrices occurring in the first and last 
expressions,
we see that $\gamma(x)=xy$, and the proof is complete.
\end{proof}

We state the corresponding result for maps into the right 
multiplier algebra
without proof.
\begin{thm}Let $X$ be an operator space and let $\psi: X\to 
I_{22}(X)$ be a
linear map. There exists $y\in I(X)^*$ such that $\psi(x)=yx$ 
for every $x\in
X$ if and only if the map $\alpha: C_2(X)\to I(S_X)$ defined by
$\alpha(\begin{pmatrix}x_1\\x_2\end{pmatrix}):=\begin{pmatrix}0 
& x_1\\0 &
\psi(x_2)\end{pmatrix}$ is completely contractive.
\end{thm}
The above theorem yields a method for determining whether or not 
a bilinear map
$m: X\times X\to X$ is a contractive quasimultiplier, i.e., 
whether or not
$m\in OAP(X)$.
\begin{cor}\label{co: bootstrap}Let $X$ be an operator space and 
let $m:
X\times X\to X$ be a bilinear map. Then the following are 
equivalent:
\begin{enumerate}
\item $m\in OAP(X)$;
\item there exists a linear map $\gamma: X\to M_l(X)$ such that 
$m(x,
y)=\gamma(x)y$ and the map $\beta: R_2(X)\to I(S_X)$ defined by 
$\beta((x_1,
x_2))=\begin{pmatrix}\gamma(x_1) & x_2\\0 & 0\end{pmatrix}$ is 
completely
contractive;
\item there exists a linear map $\psi: X\to M_r(X)$ such that 
$m(x_1, x_2)=x_1
\psi(x_2)$ and the map $\alpha: C_2(X)\to I(S_X)$ defined by
$\alpha(\begin{pmatrix}x_1\\x_2\end{pmatrix})=\begin{pmatrix}0 & 
x_1\\0 &
\psi(x_2)\end{pmatrix}$ is completely contractive.
\end{enumerate}
\end{cor}
Recall that by the results of \cite{BEZ} (see also \cite{Pa}) to 
determine
whether or not $x \to m(x_1,x)$ is a contractive left 
multiplier it is necessary and
sufficient that the map $\tau_{x_1}: C_2(X)\to C_2(X)$ defined by
$\tau_{x_1}(\begin{pmatrix}x_2\\x_3\end{pmatrix}):=\begin{pmatrix}m(x_1,x_2)\\
x_3\end{pmatrix}$ be a complete contraction. There is a similar 
result
involving $R_2(X)$ for determining when a map is a contractive 
right
multiplier.

Thus, by combining Corollary~\ref{co: bootstrap} with the 
characterizations of
multipliers, one obtains a {\bf bootstrap method} for 
determining whether or
not a bilinear map $m: X\times X\to X$ is a contractive 
quasimultiplier, i.e.,
whether or not $m\in OAP(X)$.
\begin{remark}Since we have obtained necessary and sufficient 
conditions for a
bilinear map to be in $OAP(X)$, we have in some sense given a 
full
generalization of the Blecher-Ruan-Sinclair characterization of 
unital operator
algebras \cite{BRS} to arbitrary operator algebras. However, to 
prove the
original BRS theorem by applying Corollary~\ref{co: bootstrap}, 
one still needs
to use the theory of multipliers and the proof that one obtains 
in this fashion
is not really different from the proof given in \cite{Pa}.

In the thesis of the first author \cite{Ka}, a new direct 
characterization of
the bilinear maps in $OAP(X)$ is given that is independent of 
the theory of
multipliers and is sufficiently simple that the BRS theorem can 
be deduced
directly from this characterization.
\end{remark}
The space $QMB(X)$ currently is endowed with two generally 
different norms,
$\|\cdot\|_{cb}$ and $\|\cdot\|_{qm}$. The first norm comes from 
its natural
inclusion into the space of completely bounded bilinear maps 
from $X$ into $X$
and the second from its identification with the space $QM(X)$. 
The next results
allow us to prove that for a subspace of $QM(X)$, that is 
related to ternary
structures on $X$, these two norms are the same.
\begin{thm}Let $X$ be an operator space and let 
$Y=\begin{pmatrix}X\\M_r(X)
\end{pmatrix}$ where we give $Y$ the operator space structure 
that comes from
its identification as a subspace of $I(S_X)$. Then $I(S_Y)$ can 
be identified
with the injective $C^*$-algebra $\cl 
B:=\begin{pmatrix}I_{11}(X) & I(X) &
I(X)\\I(X)^* & I_{22}(X) & I_{22}(X)\\I(X)^* & I_{22}(X) & 
I_{22}(X)
\end{pmatrix}$ in such a way that $I_{11}(Y)=I(S_X)$, 
$I(Y)=\begin{pmatrix}I(X)
\\I_{22}(X)\end{pmatrix}$ and $I_{22}(Y)=I_{22}(X)$.
\end{thm}
\begin{proof}We identify $S_Y=\begin{pmatrix}\bb C & Y\\Y^* & 
\bb C
\end{pmatrix}$ with an operator system in $\cl B$ via the map 
that sends
$\begin{pmatrix}\alpha & y_1\\y_2^* & \beta\end{pmatrix}$ where
$y_i=\begin{pmatrix}x_i\\t_i\end{pmatrix}$ to 
$\begin{pmatrix}\alpha1_{11} &
0 & x_1\\0 & \alpha1_{22} & t_1\\x_2^* & t_2^* & 
\beta1_{22}\end{pmatrix}$.
The proof of the theorem will be complete if we can show that 
any completely
positive map $\Phi: \cl B\to\cl B$ that is the identity on $S_Y$ 
must be the
identity on $\cl B$.

To this end let $\gamma: I(S_X)\to\cl B$ be defined by 
$\gamma(\begin{pmatrix}a
& b\\c & d\end{pmatrix}):=\begin{pmatrix}a & 0 & b\\0 & 0 & 0\\c 
& 0 & d
\end{pmatrix}$ and let $\delta: \cl B\to I(S_X)$ be defined by 
$\delta((a_{i,
j})):=\begin{pmatrix}a_{11} & a_{13}\\a_{31} & 
a_{33}\end{pmatrix}$.
Since $\delta\circ\Phi\circ\gamma$ is the identity on $S_X$ by 
rigidity it must
be the identity on $I(S_X)$.

To simplify notation, we define elements of $\cl B$ by
$E_{11}:=\begin{pmatrix}1_{11} & 0 & 0\\0 & 0 & 0\\0 & 0 & 
0\end{pmatrix},
E_{23}=\begin{pmatrix}0 & 0 & 0\\0 & 0 & 1_{22}\\0 & 0 & 
0\end{pmatrix}$ and
give similar definitions to $E_{22}, E_{32}, E_{33}$. Note that 
because $I(X)$
need not contain an identity we do not attempt to define 
$E_{12}, E_{13},
E_{21}$ and $E_{31}$.

We first prove that $\Phi$ fixes the five ``matrix units" 
defined above. Note
that since $\Phi$ fixes $S_Y$ we already have that
$\Phi(E_{11}+E_{22})=E_{11}+E_{22}, \Phi(E_{33})=E_{33}$ and
$\Phi(E_{23})=E_{23}$.

Since 
$\delta\circ\Phi\circ\gamma(\left(\begin{smallmatrix}1_{11} & 
0\\0 & 0
\end{smallmatrix}\right))=\left(\begin{smallmatrix}1_{11} & 0\\0 
& 0
\end{smallmatrix}\right)$, it follows that 
$\Phi(E_{11})=:P=(P_{ij})$ with
$P_{11}=1_{11}$. Since $\Phi$ is contractive and positive, it 
follows that
$P_{ij}=0$ when $(i, j)$ is (1, 2), (1, 3), (2, 1), or (3, 1) 
and that
$P_{33}\ge0$. But since $\Phi(E_{33})=E_{33}$ and
$\|\Phi(E_{11}+E_{33})\|\le1$, we have that $P_{33}=0$. Now the 
positivity of
$P$ implies that $P_{23}=P_{32}=0$.

Since $E_{23}^*E_{23}=E_{33}=\Phi(E_{23}^*E_{23})$ and 
$\Phi(E_{23})=E_{23}$,
we have that $E_{23}$ is in the right multiplicative domain of
$\Phi$, that is $\Phi(BE_{23})=\Phi(B)E_{23}, \forall B\in\cl 
B$. If we let
$\Phi(E_{22})=:Q=(Q_{ij})$, then 
$E_{23}=\Phi(E_{22}E_{23})=QE_{23}$ and it
follows that $Q_{22}=1_{22}$. But since 
$\Phi(E_{11}+E_{22})=E_{11}+E_{22}$, it
follows that $P_{22}=0$ and $Q=E_{22}$.

Thus, we have shown that these five matrix units are fixed by 
$\Phi$ as was
claimed. Since the span of these matrix units are a 
$C^*$-subalgebra of $\cl
B$, we have that $\Phi$ must be a bimodule map over this 
$C^*$-subalgebra.

Since this subalgebra contains the diagonal matrices we see that 
there exist
maps $\phi_{ij}$ such that $\Phi((B_{ij}))=(\phi_{ij}(B_{ij}))$. 
To
prove that $\Phi$ is the identity map, it will be enough to show 
that each
$\phi_{ij}$ is the identity map on its respective domain.

Using the fact that $\delta\circ\Phi\circ\gamma$ is the identity 
on $I(S_X)$
yields that $\phi_{11}$ is the identity map on $I_{11}(X)$ and 
similarly
$\phi_{13}, \phi_{31}$ and $\phi_{33}$ are the identities on 
their
respective domains.

To see that $\phi_{12}$ is the identity on its domain, note that 
for any
$u\in I(X)$ we have that \begin{multline*}\Phi(\begin{pmatrix}0 
& u & 0\\0 & 0
& 0\\0 & 0 & 0\end{pmatrix})=\Phi(\begin{pmatrix}0 & 0 & u\\0 & 
0 & 0\\0 & 0 &
0\end{pmatrix}E_{32})=\Phi(\begin{pmatrix}0 & 0 & u\\0 & 0 & 
0\\0 & 0 & 0
\end{pmatrix})E_{32}=\\\begin{pmatrix}0 & 0 & u\\0 & 0 & 0\\0 & 
0 & 0
\end{pmatrix}E_{32}=\begin{pmatrix}0 & u & 0\\0 & 0 & 0\\0 & 0 & 
0
\end{pmatrix}.\end{multline*}

Note that what was used in this argument was the bimodularity 
property of the
matrix units and the fact that certain maps were the identity 
maps. A similar
argument shows that $\phi_{21}$, $\phi_{23}$ and $\phi_{32}$ are 
all the
identity maps on their respective domains. Finally, that 
$\phi_{22}$ is the
identity follows from the rigidity of the upper left corner 
$I(\cl S_X)$ of
$\cl B$. This completes the proof of the theorem.
\end{proof} 
Given any operator space $X$ the sets $M_l(X)\cap M_l(X)^*$ and 
$M_r(X)\cap
M_r(X)^*$ are $C^*$-subalgebras of $I_{11}(Y)$ and $I_{22}(X),$
respectively. In \cite{BP} these sets were denoted $IM_l^*(X)$ 
and $IM_r^*(X)$,
respectively. They were shown to be equal to the sets $\cl 
A_l(X)$ and $\cl
A_r(X)$ of adjointable left and right multipliers, respectively, 
introduced in
\cite{B2}. We shall use the latter notation for these sets.
\begin{defn}Given an operator space $X$, we set $TER(X):=X\cap 
QM(X)^*$ and we
call this the {\bf ternary subspace of X}.
\end{defn}

Note that in the multiplication inherited from $I(S_X)$ we have 
that $TER(X) \cdot TER(X)^* \cdot TER(X) \subseteq TER(X).$
The following results give further properties of this subspace.
 
\begin{cor}Let $X$ be an operator space, let 
$Y:=\begin{pmatrix}X\\M_r(X)
\end{pmatrix}$ where we give $Y$ the operator space structure 
that comes from
its identification as a subspace of $I(S_X)$ and let $I(S_Y)$ be 
identified
with the $C^*$-algebra $\cl B$ as above. Then 
$M_l(Y)=\begin{pmatrix}M_l(X) & X
\\QM(X) & M_r(X)\end{pmatrix}$ and $\cl 
A_l(Y)=\begin{pmatrix}\cl A_l(X) &
TER(X)\\TER(X)^* & \cl A_r(X)\end{pmatrix}$.
\end{cor}
\begin{proof}One simply checks that $M_l(X)$ is exactly the 
matrix in
$I_{11}(Y)=I(S_X)$ that leaves $Y$ invariant under left 
multiplication.
\end{proof}
\begin{cor}Let $X$ be an operator space, let $x\in TER(X)$, set 
$z=x^*\in
QM(X)$ and let $m_z: X\times X\to X$ be the associated bilinear 
map, then
$\|z\|=\|m_z\|_{qm}=\|m_z\|_{cb}$.
\end{cor}
\begin{proof}We have that $\|z\|=\|m_z\|_{qm}$ by definition and 
clearly,
$\|m_z\|_{cb}\le\|z\|$. Clearly, $\|z\|=\|\begin{pmatrix}0 & 
0\\z & 0
\end{pmatrix}\|$ and this latter matrix is in $\cl A_l(Y)$. By 
Theorem~1.9 (i)
of \cite{BP} the norm of this latter matrix is equal to the norm 
of the map it
induces acting by left multiplication on $Y$, regarded as an 
element of
$CB_l(Y)$.

Hence given any $\epsilon>0$, there exists a matrix over $X, 
\|(x_{ij})\|\le1$
such that $\|(zx_{ij})\|\ge\|z\|-\epsilon$. The matrix 
$(zx_{ij})$ is a matrix
of right multipliers, and so applying \cite{BP} again, we can
find another matrix over $X$, $\|(y_{ij})\|\le1$ such that
$$\|(\sum_km_z(y_{ik}, x_{kj}))\|=\|(\sum_ky_{ik}zx_{kj})\|
\ge\|(zx_{ij})\|-\epsilon\ge\|z\|-2\epsilon.$$

Thus, we have that $\|m_z\|_{cb}\ge\|z\|-2\epsilon$ and since 
$\epsilon$ was
arbitrary, the result follows.
\end{proof}
\section{Quasicentralizers and Quasihomomorphisms}We introduce a 
new family of
bilinear maps that we call the {\em quasicentralizers} of an 
operator space and
a set of maps that we call the {\em quasihomomorphisms} and 
explore the
relationships between these maps and the space $QMB(X)$.
\begin{defn}Let $X$ be an operator space and let $m: X\times 
X\to X$ be a
bilinear map. We call $m$ a {\bf quasicentralizer}\footnote{This 
is a
generalization of a {\em quasi-centralizer} defined for an 
operator algebra
with a two-sided contractive approximate identity in \cite{Ka} 
Definition
3.2.1.} provided that there exists completely bounded maps, 
$\gamma: X\to
M_l(X)$ and $\psi: X\to M_r(X)$ such that 
$m(x,y)=\gamma(x)y=x\psi(y)$ for every
$x, y\in X$. We let $QC(X)$ denote
the set of quasicentralizers. We call a linear map $\gamma: X\to 
M_l(X)$
(respectively, $\psi: X\to M_r(X)$) a {\bf left (right)
quasihomomorphism}\footnote{This is different from a {\em 
quasi-homomorphism}
which the first author defined in \cite{Ka} Definition~3.1.1 
(2).} provided
that $\gamma(x)\gamma(y)=\gamma(\gamma(x)y)$
(respectively, $\psi(x)\psi(y)=\psi(x\psi(y)$) for every $x, y 
\in X$.
\end{defn}
These definitions are motivated by the following observations. 
If $m=m_z\in
QMB(X)$ for some $z\in QM(X)$, then $m(x, 
y)=\pi_l(x)y=x\pi_r(y)$ where
$\pi_l(x)=xz$ and $\pi_r(y)=zy$. Thus, $QMB(X)\subseteq QC(X)$. 
Moreover, the
maps $\pi_l$ and $\pi_r$ are left and right quasihomomorphisms, 
respectively.

Note that $QC(X)$ is a linear subspace of the space of bilinear 
maps from $X$
to $X$.

We begin with a few elementary observations about the 
relationships between
these concepts.
\begin{prop}Let $X$ be an operator space. A linear map $\gamma: 
X\to M_l(X)$
(respectively, $\psi: X\to M_r(X)$) is a left (respectively, 
right)
quasihomomorphism if and only if the bilinear map $m(x, 
y)=\gamma(x)y$
(respectively, $m(x, y)=x\psi(y)$) is associative. In this case,
$\gamma$ (respectively, $\psi$) is a homomorphism of the algebra 
$(X, m)$ into
$M_l(X)$ (respectively, $M_r(X)$).
\end{prop}
\begin{proof}The proof is straightforward.
\end{proof}
We shall refer to $m$ as the {\em product associated with the
quasihomomorphism.}
\begin{prop}Let $X$ be an operator space, let $\gamma : X\to 
M_l(X)$
(respectively, $\psi: X\to M_r(X)$) be a linear map and let $m(x,
y):=\gamma(x)y$ (respectively, $m(x, y):=x\psi(y)$), then
$\|m\|_{cb}=\|\gamma\|_{cb}$ (respectively, 
$\|m\|_{cb}=\|\psi\|_{cb}$).
\end{prop}
\begin{proof}By \cite{BP}, we have that the norm of a left 
(respectively,
right) multiplier is given by the cb-norm of its action as a left
(respectively, right) multiplication. Thus, given 
$\|(x_{ij})\|\le1$ and
$\|(y_{ij})\|\le1$, we have that$$\|(\sum_k m(x_{ik},
y_{kj})\|=\|(\sum_k\gamma(x_{ik})y_{kj})\|
\le\|(\gamma(x_{ij})\|\cdot\|(y_{ij})\|\le\|\gamma\|_{cb}$$and 
it follows that
$\|m\|_{cb}\le\|\gamma\|_{cb}$.

The other inequalities follow similarly.
\end{proof}
By the above results, the product associated with a completely 
contractive
quasihomomorphism is completely contractive, i.e., is in 
$CCP(X)$.
\begin{ex}This is an example of a product associated with a 
completely
contractive quasihomomorphism that is in $CCP(X)$ but not in 
$OAP(X)$.

Recall the product $m_1$ on $C_2$ of Example~\ref{ex: nonconvex} 
that is not in
$OAP(C_2)$. We have that $M_l(C_2)=M_2$ and that $\gamma: C_2\to 
M_2$ defined
by $\gamma(\begin{pmatrix}a\\b\end{pmatrix}):=\begin{pmatrix}a & 
0\\0 & b
\end{pmatrix}$ is a completely contractive quasihomomorphism 
with $m_1$ the
associated product. Thus, $\gamma$ is a completely contractive 
homomorphism of
$(C_2, m_1)$ into an operator algebra, but since $m_1$ is not an 
operator
algebra product on $C_2$, there can be no completely isometric 
homomorphism of
$(C_2, m_1)$ into an operator algebra. It is interesting to note 
that $m_1$ is
also not a quasicentralizer. In fact, it is not hard to show by 
a direct
calculation that $QC(C_2)=QMB(C_2)$.
\end{ex}
\begin{remark}Is every quasicentralizer, automatically an 
associative bilinear
map? Is every associative quasicentralizer in QMB(X)? We 
conjecture that the
answer is no to both of these questions, but we do not know of 
an example.
\end{remark}


\begin{thebibliography}{99}
\bibitem{Bl} D. P. Blecher, {\em A completely bounded 
characterization of
operator algebras}, Mathematische Annalen {\bf 303} (1995), 
227-239.
\bibitem{B2} D. P. Blecher, {\em The Shilov boundary of an 
operator space and
the characterization theorems}, Journal of Functional Analysis 
{\bf 182}
(2001), 280-343.
\bibitem{B3} D. P. Blecher, {\em One-sided ideals and 
approximate identities in
operator algebras}, Journal of the Australian Mathematical 
Society, to appear.
\bibitem{BEZ} D. P. Blecher, E. G. Effros, and V. Zarikian, {\em 
One-sided
M-ideals and multipliers in operator spaces, I}, Pacific Journal 
of Mathematics
{\bf 206 (2)} (2002), 287-319.
\bibitem{BK} D. P. Blecher and M. Kaneda, {\em The ideal 
envelope of an
operator algebra}, Proceedings of the American Mathematical 
Society, to appear.
\bibitem{BP1} D. P. Blecher and V. I. Paulsen, {\em Tensor 
products of operator
spaces}, Journal of Functional Analysis {\bf 99} (1991), 262-292.
\bibitem{BP} D. P. Blecher and V. I. Paulsen, {\em Multipliers 
of operator
spaces, and the injective envelope}, Pacific Journal of 
Mathematics {\bf 200
(1)} (2001), 1-17.
\bibitem{BRS} D. P. Blecher, Z.-J. Ruan and A. M. Sinclair, {\em 
A
characterization of operator algebras}, Journal of Functional 
Analysis {\bf
89} (1990), 188-201.
\bibitem{Br} L. G. Brown, {\em Close hereditary C*-subalgebras 
and the
structure of quasi-multipliers}, MSRI preprint, August 1985.
\bibitem{BMS} L. G. Brown, J. Mingo and N.-T. Shen, {\em 
Quasi-multipliers and
embeddings of Hilbert C*-bimodules}, Canadian Journal of 
Mathematics {\bf 46
(6)} (1994), 1150-1174.
\bibitem{CES} E. Christensen, E. G. Effros and A. M. Sinclair, 
{\em Completely
bounded multilinear maps and C*-algebraic cohomology}, 
Inventiones Mathematicae
{\bf 90} (1987), 279-296.
\bibitem{FP} M. Frank and V. I. Paulsen, {\em Injective 
envelopes of
C*-algebras as operator modules}, preprint, April 1999, Pacific 
Journal of
Mathematics, to appear.
\bibitem{H0} M. Hamana, {\em Injective envelopes of Banach 
modules}, T\^{o}hoku
Mathematical Journal (2) {\bf 30} (1978), 439-453.
\bibitem{H1} M. Hamana, {\em Injective envelopes of 
C*-algebras}, Journal of
Mathematical Society of Japan {\bf 31 (1)}, (1979), 181-197.
\bibitem{H2} M. Hamana, {\em Injective envelopes of operator 
systems},
Publications of the Research Institute for Mathematical 
Sciences, Kyoto
University {\bf 15} (1979), 773-785.
\bibitem{H3} M. Hamana, {\em Injective envelopes of dynamical 
systems},
preprint, April 1991.
\bibitem{Ka} M. Kaneda, {\em Multipliers and Algebrizations of 
Operator
Spaces}, Ph.D. Thesis, University of Houston, August 2003.
\bibitem{Pa} V. I. Paulsen, {\em Completely Bounded Maps and 
Operator
Algebras}, Cambridge Studies in Advanced Mathematics, Vol. 78, 
Cambridge
University Press, 2002.
\bibitem{Pe} G. K. Pedersen, {\em C*-algebras and Their 
Automorphism Groups},
L.M.S. Monographs, Academic Press, 1979.
\end{thebibliography}
\end{document}